\documentclass[a4paper, 11pt]{amsart}
\usepackage{graphicx} 
\usepackage{booktabs}
\usepackage{float}
\usepackage{soul}
\usepackage[citecolor=Navy]{hyperref}
\usepackage[ruled,vlined]{algorithm2e}
\usepackage{amsfonts,amssymb,amsmath,amsthm}
\usepackage{algorithmic}
\usepackage{amssymb,amscd,amsxtra,calc}
\usepackage{cmmib57}
\usepackage{graphicx}
\usepackage{mathrsfs,enumerate,xcolor}

\usepackage[all]{xy}
\usepackage{multirow} 
\usepackage{psfrag,xmpmulti,amscd,color,pstricks, import}
\usepackage{tikz-cd} 

\newtheorem{thm}{Theorem}[section]
\newtheorem{cor}[thm]{Corollary}

\theoremstyle{definition}
\newtheorem{defn}[thm]{Definition}

\newtheorem*{rem*}{Remark}

\title[Riemann Hypothesis and dynamics of BNQN]{The Riemann hypothesis and dynamics of Backtracking New Q-Newton's method}
\date{\today}

\author[Tran]{Thuan Quang Tran}
\address{Department of Mathematics, University of Oslo, Norway}
\email{thuanqt@math.uio.no}

\author[Truong]{Tuyen Trung Truong}
\address{Department of Mathematics, University of Oslo, Norway}
\email{tuyentt@math.uio.no}

\begin{document}

\maketitle

\begin{abstract}

A new variant of Newton's method - named Backtracking New Q-Newton's method (BNQN) -  was recently introduced by the second author. This method has good convergence guarantees, specially concerning finding roots of meromorphic functions. This paper explores using BNQN for the Riemann xi function. We show in particular that the Riemann hypothesis is equivalent to that all attractors of BNQN lie on the critical line. We also explain how an apparent relation between the basins of attraction of BNQN and Voronoi's diagram can be helpful for verifying the Riemann hypothesis or finding a counterexample to it. Some illustrating experimental results are included, which convey some interesting phenomena. The experiments show that BNQN works very stably with highly transcendental functions like the Riemann xi function and its derivatives. Based on insights from the experiments, we discuss some concrete steps on using BNQN towards the Riemann hypothesis, by combining with de Bruijn -Newman's constant. Ideas and results from this paper can be extended to other zeta functions.  
\end{abstract}

\section{Introduction}

The Riemann hypothesis is a famous long standing open problem in mathematics. It concerns the zeros of the Riemann zeta function $\zeta (s)$, which is defined in the domain $\mathcal{R}(s)>1$ (where $\mathcal{R}(.)$ is the real part of a complex number) by an absolutely convergent power series: 

\begin{eqnarray*}
\zeta (s) =\sum _{n=1}^{\infty}\frac{1}{n^s},
\end{eqnarray*}
where $n^s=e^{s \log (n)}$. The Riemann zeta function has an analytic continuation for \(\mathcal{R}(s) > 0\) by presenting it as a Riemann--Stietjes integral \cite{RefMontVaug}: 
\begin{equation}
\zeta(s) =  \frac{s}{s-1} +s\int_1^{\infty} \{u\}u^{-s-1}\;du, 
\end{equation}
where $\{u\}$ is the fractional part of \(u\) (for example, if \(u = 2.52\), then \(\{2.52\} = 0.52\)).

The integral formula helps to establish that the Riemann zeta function is a meromorphic function with only one simple pole at $s=1$ with residue \(1\). Moreover, one has the following functional equation: 

\[
    \zeta(s) = 2^s\pi^{s-1}\sin\left(\frac{\pi s}{2}\right)\Gamma(1-s)\zeta(1-s), \quad s \in \mathbb{C}\setminus\{1\}.
\]

Some zeros of the Riemann hypothesis are easy to determine, and hence named trivial zeros. They are precisely those of the form: $-2n$, $n\in \mathbf{N}_{>0}$, which are related to the Gamma function. The other zeros of the Riemann zeta function are called non-trivial zeros. 

{\bf The Riemann hypothesis}: Non-trivial zeros of the Riemann zeta function all lie on the critical line $\mathcal{R}(s)=0.5$. 

For the purpose of this paper, it is more convenient to work with the Riemann xi function: 

\[
    \xi(s) = \frac{s(s-1)}{2}\pi^{-\frac{s}{2}}\Gamma\left(\frac{s}{2}\right)\zeta(s), \quad s \in \mathbb{C}.
\]

The Riemann xi function is an entire function of order 1, and satisfies a symmetric relation (functional equation): $\xi (s)=\xi (1-s)$. This relation helps to show that $\xi$ has real values on the critical line. 

The roots of the Riemann xi function are precisely the non-trivial roots of the Riemann zeta function. Hence, the Riemann hypothesis is then the statement that all zeros of the Riemann xi function belong to the critical line. 

There have been a lot of information known about the Riemann hypothesis. For example, by Euler's product it is known that zeros of the Riemann zeta function satisfy $\mathcal{R}(s)\leq 1$, and a more complicated argument shows that $\mathcal{R}(s)<1$ (Hadamard and de la Vall\' ee-Poussin, independently, see \cite{RefMontVaug}). Hence, by the functional equation, non-trivial roots of the Riemann zeta function lies in the critical strip $0<\mathcal{R}(s)<1$.  Moreover, a positive portion of the non-trivial zeros of the Riemann zeta function is shown to be on the critical line (by work of Hardy and Littlewood, Selberg, Levinson and others, the most current records are \cite{RefCon} and \cite{RefPratt}, where respectively at least $2/5$ and $5/12$ of the zeros are shown to be on the critical line). Bohr and Landau showed that for any $\epsilon >0$ there are at most $O(T)$ roots of the Riemann zeta function with $\mathcal{R}(s)>0.5+\epsilon $ and $|\mathcal{I}(s)|<T$ (negligible to the number of the zeros of the Riemann zeta function in the domain $0<\mathcal{R}(s)<1$ and $|\mathcal{I}(s)|<T$), see \cite{RefIvic} for improvements on this. Additionally, it is known that a positive portion of the non-trivial roots of the Riemann zeta function are simple (i.e. of multiplicity 1), the current record is $41\%$  \cite{RefBui}. Moreover, there are works (in particular by Levison and Conrey \cite{RefLevinson}\cite{RefCon2}) which show that the density of zeros of the m-th derivative of the Riemann $\xi$ function on the critical line tends to $1$ when $m\rightarrow \infty$.  Also, it is verified by numerical methods that the Riemann hypothesis is true for zeros $s$ with large imaginary parts (starting with works by Riemann himself, Turing's method for finding roots on the critical line, and specially works with very large heights by Odlyzko \cite{RefOdlyzko}), the current record being  $|s|\leq 3\times 10^{12}$ \cite{RefPlattTrudgian} (where a good reference list on this topic can be found). 

The Riemann hypothesis has  wide connections to many fields in mathematics and physics. As such, there are many approaches and equivalences to it (see e.g. \cite{RefBroughan1}, \cite{RefBroughan2},\cite{RefBroughan3}). This paper concerns the connections between the Riemann hypothesis and root finding algorithms, specifically Newton-type methods. We note that this direction is different from the works mentioned above on verifying the Riemann hypothesis up to a given height: those works are based on Turing's method, and concentrate on the critical line $\mathcal{R}(s)=0.5$ only. The Newton-type's methods, on the other hand, work on the whole critical strip $0\leq \mathcal{R}(s)\leq 1$, and can potentially produce explicit counterexamples in case the Riemann hypothesis is not true. According to Wikipedia, while the numerical evidence for Riemann hypothesis is abundant, analytic number theory has many conjectures (e.g. Skewes number) which are true up to a very large number, but are actually false. It is known that there are zeta functions, e.g. Davenport-Heilbroon function \cite{davenport-heilbronn}\cite{spira}\cite{balanzario-ortiz}, for which the Riemann hypothesis fails.  

We recall that Newton's method, for finding roots of a meromorphic function $g:\mathbf{C}\rightarrow \mathbf{C}\cup \{\infty\}$, is the following iterative method: 
\begin{eqnarray*}
z_{n+1}=z_n-\frac{g(z_n)}{g'(z_n)}, 
\end{eqnarray*}
starting from an initial point $z_0\in \mathbf{C}$. One hopes that $\lim _{n\rightarrow\infty}z_n$ exists and is a root of $g(z)$. However, it is known that even for polynomials of small degrees, this is not always possible, e.g. the fractal nature of the basins of attraction (please consult the Wikipedia page on Newton fractal for details). 

Associated to it is Newton's flow: 
\begin{eqnarray*}
z'(t)=-g(z(t))/g'(z(t)),
\end{eqnarray*}
with initial value $z(0)=z_0$. One hopes that $\lim _{t\rightarrow\infty}z(t)$ exists and is a root of $g(z)$. Indeed, strong convergence results are known in case $g(z)$ is a rational function (i.e. the quotient of two polynomials) \cite{RefMe}\cite{RefHK}, but not much is known in case $g(z)$ is transcendental (in particular, for the Riemann zeta/xi functions). Another issue with Newton's flow is that one cannot really use it in practice, but needs to use some approximation schemes (like Euler's scheme, for which Newton's flow becomes the Relaxed Newton's method below), and hence results proven for Newton's flow may not transfer to its approximations. 

A direct variant of Newton's method is Relaxed Newton's method 
\begin{eqnarray*}
z_{n+1}=z_n-\alpha g(z_n)/g'(z_n),
\end{eqnarray*}
and its randomized version:
\begin{eqnarray*}
z_{n+1}=z_n-\alpha _ng(z_n)/g'(z_n),
\end{eqnarray*}
starting from an initial point $z_0$ (where $\alpha _n$ are complex numbers randomly chosen in an appropriate domain). One hopes that $\lim _{n\rightarrow\infty}z_n$ converges and is a root of $g(z)$. If $g(z)$ is a generic rational function and $0\not= \alpha$ is small enough (depending on $g$), then \cite{RefHK} shows that Relaxed Newton's method has strong convergence guarantee (however, the convergence may be very slow). If $g(z)$ is a polynomial, then \cite{RefS} shows that Random Relaxed Newton's method, with a suitable randomization of the values $\alpha _n$, has convergence guarantee. Again, not much is known about the convergence guarantees of these methods for transcendental functions, in particular for the Riemann zeta/xi functions.  

Prior to this paper, there have been experimental works applying Newton's method and Newton's flow to finding roots of the Riemann zeta function, see for example \cite{RefNeuberger2}\cite{RefNeuberger1}. There have been also experiments applying Newton's method to the Riemann xi function \cite{RefKawahira}. There, it is seen that in large scale Newton's method and Relaxed Newton's method for the Riemann xi function looks quite simple, however it is very fractal in the critical strip. Some experiments in this paper confirm this fractal feature of Newton's method. We note that there are many  other variants of Newton's method, in both equation solving and optimization, such as Levenberg-Marquardt method and Regularized Newton's method, but the existing literature concerns mainly local convergence properties, while for global convergence guarantees they require very strong conditions (e.g. convextity or Lipschitz continuity of the gradient or Hessian matrix) which are not satisfied by the Riemann zeta or xi function. We are not aware of any work applying these  methods to the Riemann hypothesis.  

Note that Newton's method and Newton's flow should behave well when applied to the Riemann xi function. Indeed, by the result in \cite{RefHinkkanen}, the Riemann hypothesis is equivalent to that $\mathcal{R}(\xi (z)/\xi '(z))<0$ for $\mathcal{R}(z)<0.5$. Therefore, if one starts from an initial point of the form $z_0=(x_0,y_0)$ where $x_0<0.5$ and apply Newton's method, then if the constructed sequence $z_n=(x_n,y_n)$ satisfies $x_n<0.5 $ for all $n$, and the Riemann hypothesis holds, then $x_{n+1}>x_n$ for all $n$. A similar claim for Newton's flow. There are some informal discussions on the dynamics of Relaxed Newton's method for the Riemann xi function (as well as various other dynamical systems) in \cite{RefHinkkanen}. \cite{RefSchr} provides a sketch of idea on how to show that Newton's method applied to the Riemann xi function with initial points of the form $6+i cn/log n$, where $c$ is a constant and $n\in \mathbf{Z}\backslash \{0\}$ can find all roots of the $\xi$ function. However, global convergence is  unknown, for example it is unclear whether Newton's method applied to the Riemann xi function could have some strange attractors (e.g. attracting periodic cycles, as seen in the case of polynomials). 

For a meromorphic function $g(z)$, \cite{RefKawahira}  defines $\nu _g(z)=z - \frac{g(z)}{zg'(z)}$. Concerning theoretical connections between the Riemann hypothesis and dynamical systems, there is the following interesting result. 
\begin{thm}[{\cite{RefKawahira}}] 

A) The Riemann hypothesis is equivalent to that there is no topological disk $D$ contained in the critical strip such that $\overline{N_{\zeta}(D)}\subset D$, where $N_g(z)=z-g(z)/g'(z)$ is Newton's method for $g$. 
 
B)  The following  statements are equivalent (and similarly in the case we replace $\zeta$ by $\xi$): 

Statement 1: The Riemann hypothesis is true and the non-trivial zeros are simple.

Statement 2: The meromorphic function $\nu _{\zeta }(z)$ has no attracting fixed point. 

\label{TheoremKawahira}
\end{thm}

By Statement 2 in part B of the above theorem, one sees that $\nu _{g}$ is not a good root finding method, specially for $\zeta$ or $\xi$. 

Another result which uses the Riemann xi function, indirectly related to Newton's flow, is: 

\begin{thm}[{\cite{RefSch}}]  The following statements are equivalent. 

Statement 1: The Riemann hypothesis holds. Moreover, all zeros of the Riemann xi function and of its first derivative are simple. 

Statement 2: All lines of constant phase of the Riemann xi function corresponding to $\pm \pi$, $\pm 2\pi$, $\pm 3\pi $, $\ldots$ merge with the critical line. 

\label{TheoremSch}\end{thm}

The main focus of this paper is on a relatively new variant of Newton's method called Backtracking New Q-Newton's method \cite{RefT}. It is indeed a variant of Newton's method for optimization, here when applied to finding roots of a meromorphic function $g(z)$ in 1 complex variable $z$, one reduces to finding global minima of an associated function $F(x,y)=|g(x+iy)|^2/2$ from $\mathbf{R}^2$ to $\mathbf{R}\cup \{\infty \}$. It has very strong convergence guarantee for finding roots of meromorphic functions, and hence very appropriate to use to find roots of the Riemann zeta/xi functions or other zeta functions. Please see Section 2 for the precise algorithm and for some main properties of this method. 

Let $\phi :\mathbf{C}\rightarrow \mathbf{C}$ be a function (not necessarily continuous). We use the following simple notion of attractors: 

{\bf Attractor:} An attractor of $\phi$ is a non-empty compact subset $S$ of $\mathbf{C}$ which satisfies $\phi (S)=S$ and which has an open neighbourhood $S\subset U$ so that for all $x\in U$, $\lim _{n\rightarrow \infty}dist(\phi ^n(x),S)=0$. To ease the discussion, we also require that $S$ is minimal, in the sense that it does not contain a smaller non-empty subset which also satisfies the condition in the previous sentence. 

The basin of attraction for an attractor $S$ consists of initial points $z_0\in \mathbf{C}$ whose orbit under the dynamics of $\phi$ converges to $S$.  

Our main theoretical results are the following two theorems.  

\begin{thm}
The following 3 statements are equivalent: 

Statement 1: The Riemann hypothesis is true. 

Statement 2: All attractors of the dynamics of BNQN applied to the Riemann xi function are contained in the critical line. 

Statement 3: For any initial point $z_0$, the constructed sequence $\{z_n\}$ will either converge to a point on the critical line or to the point $\infty$.

\label{TheoremMain}\end{thm}

Under stronger assumptions, we have the following result, which shows that BNQN is a good root finding method for $\xi$. 

\begin{thm} Assume that the Riemann hypothesis holds, and moreover all roots of the Riemann xi function and its derivative are simple.  If the parameters of BNQN are randomly chosen, then there is an exceptional set $\mathcal{E}\subset \mathbf{C}$ of Lebesgue measure 0 so that the following is true. If $z_0\in \mathbf{C}\backslash \mathcal{E}$ is an initial point, then for BNQN applied to the Riemann xi function, the constructed sequence $\{z_n\}$ satisfies the following two alternatives: 

Alternative 1: $\lim _{n\rightarrow\infty}z_n=z^*$, where $z^*$ is a root of the Riemann xi function, and  the rate of convergence is quadratic. 

Alternative 2: $\lim _{n\rightarrow\infty}z_n=\infty$. 

\label{Theorem2}\end{thm} 

\begin{figure}
    \centering
    \includegraphics[width=.45\linewidth]{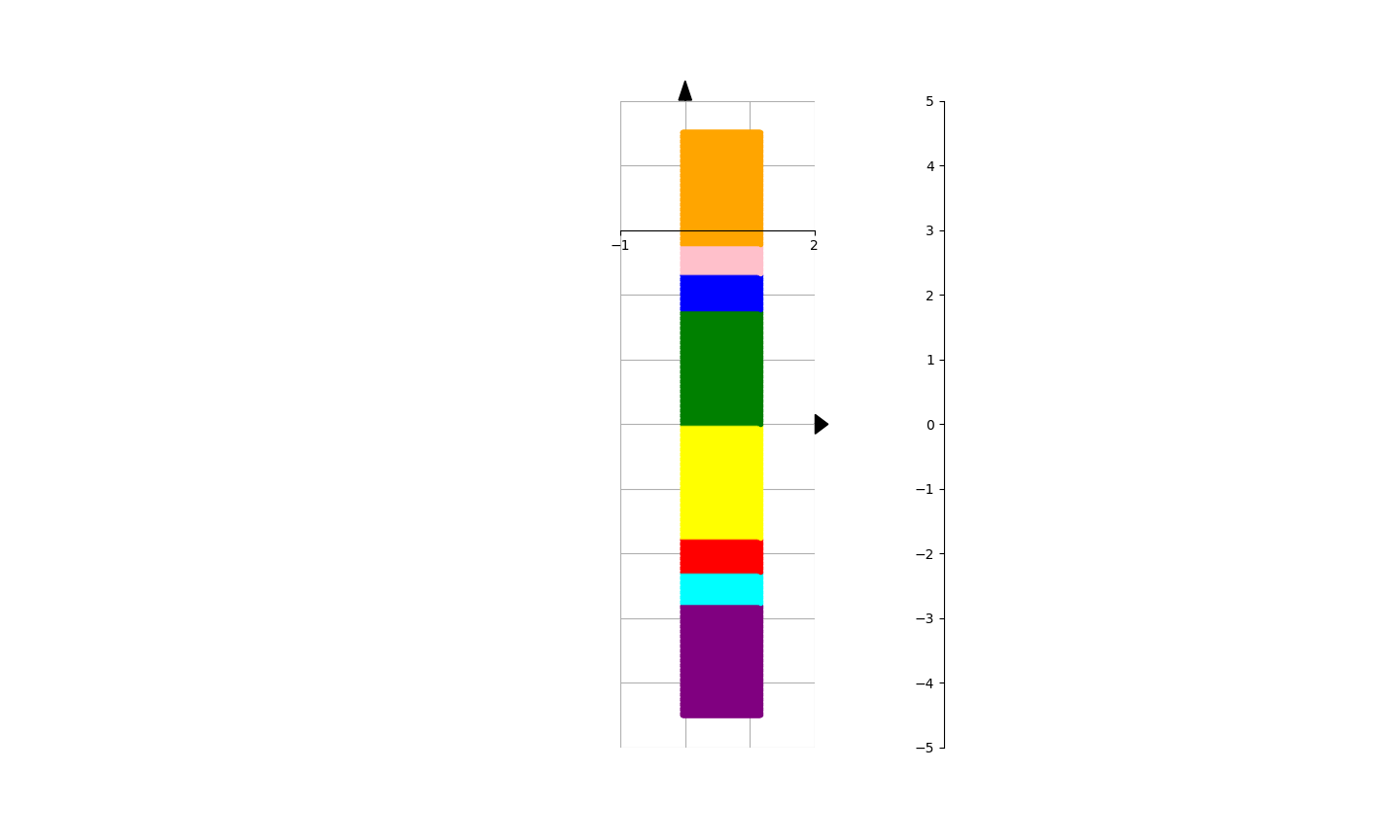}
    \includegraphics[width=.45\linewidth]{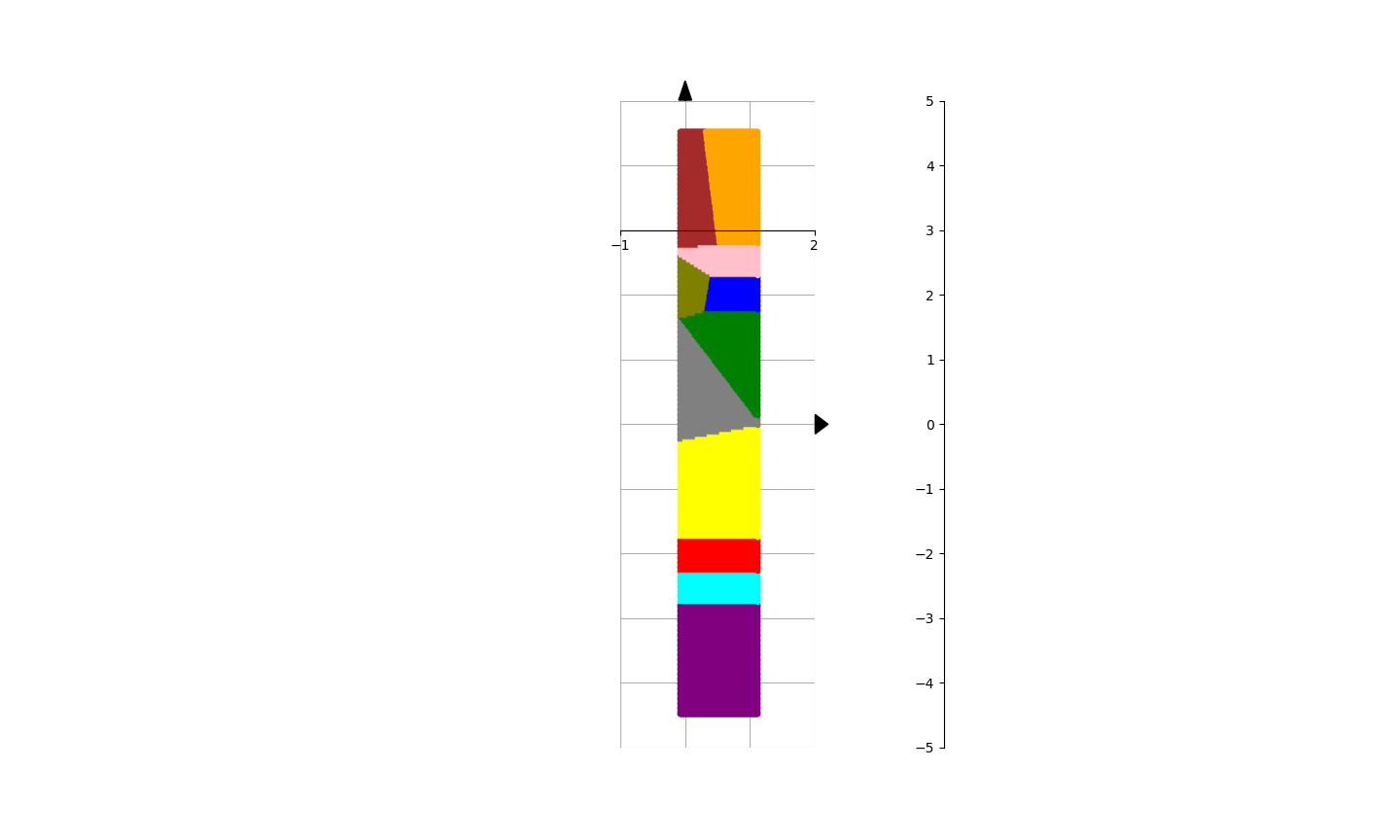}

    \bigskip
    
    \includegraphics[width=1.5cm]{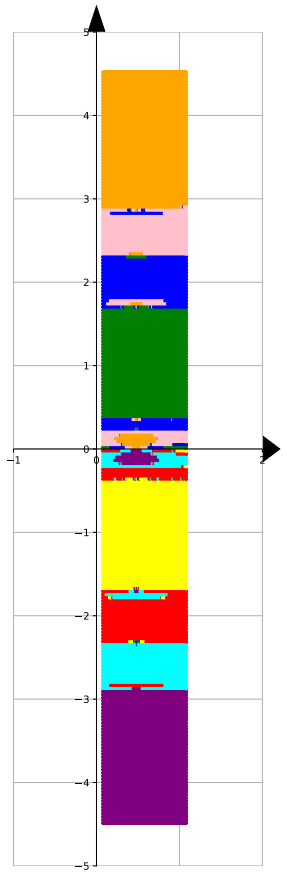}
    \includegraphics[width=1.5cm]{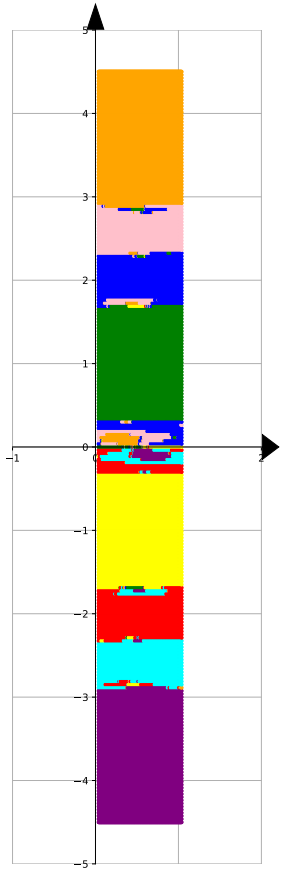}

    \bigskip
   \includegraphics[width=7cm]{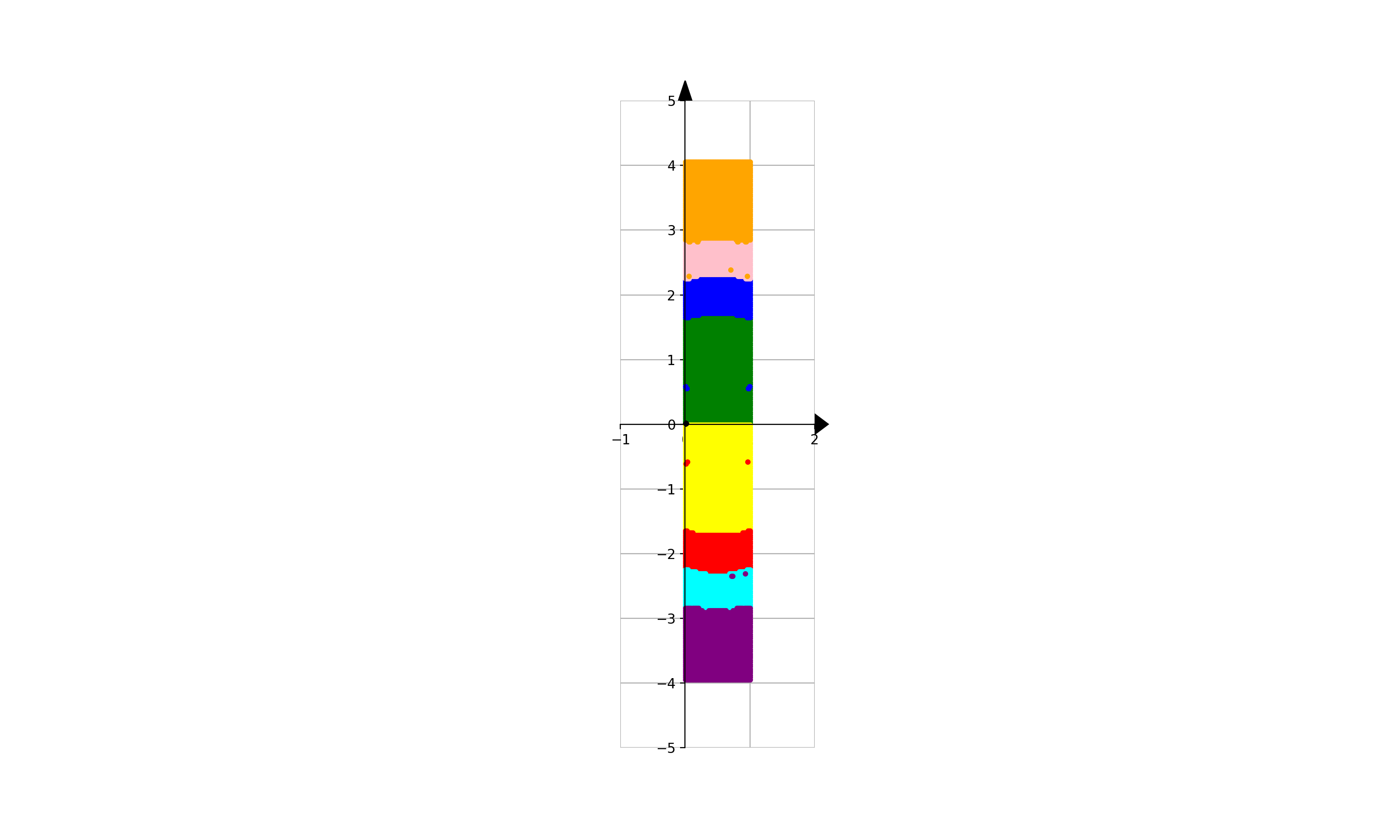}
    
    \caption{Basins of attraction for finding roots of the polynomial of degree 8, whose roots are the first 8 roots of the Riemann xi function: Root 1 $\sim$ $0.5+14.13472514173i$ (points in the basin of attraction has green colour), Root 2 $\sim$ $0.5-14.13472514173i$ (yellow), Root 3 $\sim$ $0.5+21.02203963877i$ (blue), Root 4 $\sim$ $0.5-21.02203963877i$ (red), Root 5 $\sim$ $0.5+25.01085758014i$ (pink), Root 6 $\sim$ $0.5-25.01085758014i$ (cyan), Root 7 $\sim$ $0.5+30.42487612585i$ (orange), Root 8 $\sim$ $0.5-30.42487612585i$ (purple). Other points have black colour. The y-axis is rescaled by $0.1$, hence $y=1$ in the pictures correspond to $y=10$ in reality.  Pictures are referenced to from top to bottom, from left to right. Row 1: left picture is Voronoi's diagram for the 8 roots, right picture is Voronoi's diagram for the 8 roots and some extra points which are off the critical line. Row 2: left picture is for Newton's method, right picture is for Random Relaxed Newton's method. Row 3: BNQN. One can see that pictures for Newton's method and Random Relaxed Newton's method display fractal features.}
    \label{fig:f1}
\end{figure}

Note that in \cite{fornaess-etal2}, it was observed that the basins of attraction produced by BNQN look similar to the Voronoi's diagram of the roots. Here we recall that the Voronoi diagram \cite{RefV1}\cite{RefV2} of a discrete set of points $\{z_n^*\}$ in $\mathbf{C}$ is the union of Voronoi's cells, where each Voronoi cell $V(z_n^*)$ consists of points $z\in \mathbf{C}$ for which $|z-z_n^*|<\min _{j\not= n}|z-z_j^*|$. The boundaries of Voronoi's cells are line segments. Note that if the whole sequence $\{z_n^*\}$ belongs to a line $L$, then the boundaries of the Voronoi cells are lines orthogonal to $L$ and go through the middle points of the intervals connecting consecutive points in the sequence $\{z_n^*\}$. Figure \ref{fig:f1} presents the Voronoi's diagram for the first 8 roots of the Riemann xi function in comparison to the basins of attraction found by Newton's method, Random Relaxed Newton's method and BNQN for the polynomial of degree $8$ whose roots are exactly the first 8 roots of the Riemann xi function. 

In general, this similarity between Voronoi's diagrams and basins of attraction for BNQN should be interpreted as follows. The boundary curves of the basins of attraction should go through the critical points of the meromorphic function $g(z)$, and they do not need to be literally line intervals but  are curved, depending on the flows of the gradient of $F(x,y)$ and related vector fields. For nice functions like $g(z)=\sin (z)$, where all roots are on the real axis and all critical points are precisely middle points of consecutive roots, the picture of basins of attraction is precisely that of Voronoi's diagram. On the other hand, for other functions like the one in Figure \ref{fig:f1} or a very complicated one like the Riemann xi function (see Section 3), the similarity must be viewed in a more broad sense. 

In any case,  this suggests the following idea to finding all roots of the Riemann xi function which are closest to the line $x=0$:   

{\bf Idea.} We consider initial points $(0,y_j)$, where $y_j$ is randomly chosen, and run BNQN. We should be able to find all roots of the Riemann zeta function which are closest to the line $x=0$ (in particular a counter example in case the Riemann hypothesis fails). 

The organisation of this paper is as follows. In Section 2 we recall some necessary tools, and then prove Theorems \ref{TheoremMain} and \ref{Theorem2}. In Section 3 we present some experiments concerning applied BNQN to finding roots of the Riemann xi function. Various interesting phenomena from experiments are given in Section 4, where we elaborate on the above idea to discuss some concrete steps on using BNQN towards the Riemann hypothesis, by combining with de Bruijn -Newman's constant. 

{\bf Acknowledgments}: The authors thank Juan Arias de Reyna and Tomoki Kawahira for helping with some references, {Viktor Balch Barth for translation of some articles from German}, and Terje Kvernes for helping with experiments. This paper is part of the first author's Master's thesis. The second author is partially supported by Young Research Talents grant 300814 from Research Council of Norway.

\section{Proofs of theoretical results}

First we recall the definition of BNQN, and some of its main properties which are useful for the purpose of this paper. We also recall about the Laguerre-P\'olya class. Then, we provide the proofs of the main theorems. 

\subsection{Backtracking New Q-Newton's method} Finding roots of a meromorphic function $g(z)$ can be reduced to finding global minima of an associated function as follows. Let $z=x+y$ where $x,y\in \mathbf{R}$ and define $F(x,y)=|g(x+iy)|^2$. We note that $F\geq 0$, hence if $z^*=(x^*,y^*)$ is a root of $g(z)$, then $(x^*,y^*)$ is a global minimum of $F$. Thus, one can apply an optimization method to find roots of $g(z)$. 

If one wants to find the minima of a function $H:\mathbf{R}\rightarrow\mathbf{R}$, then the first step would be to find critical points $H'(x)=0$. When one applies Newton's method to find roots of $H'$, one obtains the scheme $z_{n+1}=z_n-H"(z_n)^{-1}.H'(z_n)$. In general, the version of Newton's method for optimization (for a function $F$ as in the previous paragraph) is as follows: 
\begin{eqnarray*}
z_{n+1}=z_n-(\nabla ^2F(z_n))^{-1}.\nabla F(z_n),    
\end{eqnarray*}
starting from an initial point $z_0=(x_0,y_0)$, provided the matrix $\nabla ^2F(z_n)$ is invertible. Here, $\nabla F$ is the gradient, and $\nabla ^2F$ is the Hessian matrix. We hope that the sequence $z_n$ will converge to a root of $g(z)$. 

Like (direct) Newton's method, this optimization version also has issues with convergence guarantee. For example, it cannot avoid saddle points (this can be observed with such simple functions as $g(z)=z^2-1$, where $0$ is the saddle point of the function $F=|g|^2/2$). 

In \cite{RefTT}, a new variant of Newton's method for optimization, New Q-Newton's method (NQN), is proposed. A crucial idea is that we should change the sign of negative eigenvalues of the Hessian matrix to positive, in order to avoid saddle points. It can avoid saddle points and still has fast rate of convergence near non-degenerate local minima, but it still does not have global convergence guarantee. An improvement of NQN is proposed in \cite{RefT}, which incorporates also Armijo's Backtracking line search \cite{RefAr} (a well known technique to boost convergence for algorithms), and is named Backtracking New Q-Newton's method (BNQN).  

The following is a modification from \cite{RefFHTW}:

\medskip
{\color{blue}
 \begin{algorithm}[H]
\SetAlgoLined
\KwResult{Find a minimum of $F:\mathbb{R}^m\rightarrow \mathbb{R}$}
Given: $\{\delta_0,\delta_1,\ldots, \delta_{m}\} \subset \mathbb{R}$, $\theta \geq 0$, $0<\tau $ and $0<\gamma _0\leq 1$;\\
Initialization: $z_0\in \mathbb{R}^m$\;
$\kappa:=\frac{1}{2}\min _{i\not=j}|\delta _i-\delta _j|$;\\
 \For{$k=0,1,2\ldots$}{ 
    $j=0$\\
  \If{$\|\nabla F(z_k)\|\neq 0$}{
   \While{$minsp(\nabla^2F(z_k)+\delta_j \|\nabla F(z_k)\|^{\tau}Id)<\kappa  \|\nabla F(z_k)\|^{\tau}$}{$j=j+1$}}
  
 $A_k:=\nabla^2F(z_k)+\delta_j \|\nabla F(z_k)\|^{\tau}Id$\\
$v_k:=A_k^{-1}\nabla F(z_k)=pr_{A_k,+}(v_k)+pr_{A_k,-}(v_k)$\\
$w_k:=pr_{A_k,+}(v_k)-pr_{A_k,-}(v_k)$\\
$\widehat{w_k}:=w_k/\max \{1,\theta \|w_k\|\}$\\
$\gamma :=1$\\
 \If{$\|\nabla F(z_k)\|\neq 0$}{
   \While{$F(z_k-\gamma \widehat{w_k})-F(z_k)>-\gamma \langle\widehat{w_k},\nabla F(z_k)\rangle/3$}{$\gamma =\gamma /2$}}

$z_{k+1}:=z_k-\gamma \widehat{w_k}$
   }
  \caption{Backtracking New Q-Newton's method (BNQN)} \label{table:alg2}
\end{algorithm}
}
\medskip

Both NQN and BNQN can be applied to any dimensions and any $C^2$ function. Concerning finding roots of meromorphic functions, we have the following result \cite{RefT}\cite{RefTT}: 

\begin{thm} Let $g(z):\mathbf{C}\rightarrow \mathbf{P}^1$ be a non-constant meromorphic function. 
Define a function $F:\mathbf{R}^2\rightarrow [0,+\infty]$ by the formula $F(x,y)=|g(x+iy)|^2/2$. 

Given an initial point $z_0\in \mathbf{C}$, which is not a pole of $g$, we let $\{z_n\}$ be the sequence constructed by BNQN applied to the function $F$ with initial point $z_0$. If the objective function has compact sublevels (i.e. for all $C\in \mathbf{R}$ the set $\{(x,y)\in \mathbf{R}^2:~F(x,y)\leq C\}$ is compact), we choose $\theta \geq 0$, while in general we choose $\theta >0$.

1)  If $F$ has compact sublevels, then $\{z_n\}$ converges. 

2) If $\{z_n\}$ has a bounded subsequence, then $\{z_n\}$ converges to a point $z^*$, which is a root of $f(z)f'(z)$. 

3) Assume that the parameters $\delta _0,\delta _1,\delta _2$ in BNQN are randomly chosen. Assume also that $g(z)$ is generic, in the sense that $\{z\in \mathbf{C}:~g(z)g''(z)=g'(z)=0\}=\emptyset$. There exists an exceptional set $\mathcal{E}\subset \mathbf{C}$ of zero Lebesgue measure so that if $z_0\in \mathbf{C}\backslash \mathcal{E}$, then $\{z_n\}$ must satisfy one of the following two options: 

Option 1: $\{z_n\}$ converges to a root $z^*$ of $g(z)$, and if $\gamma _0=1$ in the algorithm then the rate of convergence is quadratic. 

Option 2: $\lim _{n\rightarrow\infty}|z_n| =+\infty$. 

\label{TheoremMeromorphic}\end{thm}

In terms of computational resources, ensuring Armijo's condition (the second While loop in the algorithm) is the most expensive aspect of BNQN. For example, if $\gamma$ (the learning rate) $=0.0625=(1/2)^4$ after the While loop, it means that we had to check Armijo's condition $4$ times. Note that by part 4 of Theorem \ref{TheoremMeromorphic}, when we are close enough to a non-degenerate root, then $\gamma =1$ after the While loop, that is we only had to check Armijo's condition once.

\subsection{Laguerre-P\'olya class} 
\cite{RefBruijn}\cite{RefNewman} defines a family of functions \(\{H_t\}_{t \in \mathbf{R}}\)
\[
    H_t(z) = \int_{0}^{\infty} \Phi(u)e^{tu^2}\cos(zu)\;du,
\]
where
\[
    \Phi(u) = \sum_{n=1}^{\infty} \left(2\pi^2n^4e^{9u} - 3\pi n^2e^{5u}\right)\exp\left(-\pi n^2e^{4u}\right).
\]
There is a constant $-\infty <\Lambda < 1/2$ (de Bruijn--Newman constant) so that the function \(H_{t}\) has only real zeros if and only if \(t \geq \Lambda\). Otherwise, the function \(H_{t}\) has some non-real zeros for \(t < \Lambda\). We have the following relation between the Riemann xi function and the function \(H_0\).
\[
    \frac{1}{8}\xi\left(\frac{1}{2} + \frac{iz}{2}\right) = H_0(z).
\]
By the relation, RH is true if and only if \(H_0\) has only real zeros. Moreover, it is known that the Riemann hypothesis holds if and only if $\Lambda = 0$ \cite{RefRodgersTao}\cite{RefPolymath}. One can also improve the upper bound $\Lambda$ in terms of a height $H$ for which the Riemann hypothesis is valid in the domain $\{z\in \mathbf{C}: ~0\leq \mathcal{R}(z)\leq 1,~|\mathcal{I}(z)|\leq H\}$ \cite[Table 1]{RefPolymath}. Therefore, the Riemann hypothesis is intimately related to the following class of entire functions. 
\begin{defn}
Let \(f \colon \mathbb{C} \to \mathbb{C}\) be an entire function. Then \(f\) belongs to the Laguerre--P\'olya class, denoted by \(\mathcal{LP}\), if \(f\) can be expressed as the Hadamard factorization
\[
    f(z) = Az^me^{-az^2 + bz}\prod_{k=1}^{\infty} \left(1 - \frac{z}{z_k}\right)e^{\frac{z}{z_k}}, \quad A,b \in \mathbb{R}\setminus\{0\},
\]
where \(a \geq 0\), \(m \in \mathbb{N}\) and \(z_k \in \mathbb{R}\setminus \{0\}\) are  zeros of \(f\) such that 
\[
    \sum_{k=1}^{\infty} \frac{1}{z_k^2} < \infty.
\]
\end{defn}
One has the following useful characterisation for that \(f \in \mathcal{LP}\).
\begin{thm}[{\cite[Laguerre--P\'olya Theorem]{RefPolya}}]
Let \(f \colon \mathbb{C} \to \mathbb{C}\) be a non-constant entire function. Then \(f \in \mathcal{LP}\) if and only if there exists a sequence of (complex) polynomials \(\{P_n\}_{n \in \mathbb{N}}\) with only real zeros such that it converges uniformly to \(f\) in  \(|z| \leq R\), for every \(R > 0\).
\end{thm}

More precisely, if $f(z) = Az^me^{-az^2 + bz}\prod_{k=1}^{\infty} \left(1 - \frac{z}{z_k}\right)e^{\frac{z}{z_k}}$, then we can approximate it by 
\begin{eqnarray*}
Az^me^{-az^2 + b_nz}\prod_{k=1}^{n}\left(1 - \frac{z}{z_k}\right),
\end{eqnarray*}
where $b_n = b + \sum_{k=1}^n \frac{1}{z_k}$. Then we use 
\begin{eqnarray*}
\lim _{k\rightarrow\infty}\left(1 - \frac{az^2}{k}\right)^{k}=e^{-az^2},
\end{eqnarray*}
and
\begin{eqnarray*}
\lim _{k\rightarrow\infty}\left(1 + \frac{b_nz}{k}\right)^{k}=e^{b_nz},
\end{eqnarray*}
to approximate $e^{-az^2+bz}$ by polynomials of the form (where $k_n$ is relatively large in comparison to $n$ and $b_n$)
\begin{eqnarray*}
\left(1 - \frac{az^2}{k_n}\right)^{k_n}\left(1 + \frac{b_nz}{k_n}\right)^{k_n}. 
\end{eqnarray*}
Since $a\geq 0$ and $b\in \mathbf{R}$ (and hence also $b_n\in \mathbf{R}$), the above polynomials have only real roots and have real values on the real axis. 

\subsection{Proofs of main theoretical results} Now we are ready to prove Theorems \ref{TheoremMain} and  \ref{Theorem2}. 

\begin{proof}[Proof of Theorem \ref{TheoremMain}] We will present the equivalence between Statement 1 and Statement 2. The other equivalences are similar. 

($\Rightarrow $) We first prove that Statement 1 implies Statement 2. Indeed, Statement 1 implies that $H_0(s)$ is of Laguerre-P\'olya class. Therefore, there is a sequence of polynomials which have only real zeros and which converge locally uniformly to $H_0(s)$. From the relations between $\xi $ and $H_0$, this implies that there is a sequence of polynomials $P_n(s)$ which have only zeros on the critical line and which converges locally uniformly to $\xi $. Then the derivatives $P_n'(s)$ converges locally uniformly to $\xi '(s)$. By Gauss-Lucas' theorem, roots of $P_n'(s)$ lie in the convex hull of roots of $P_n(s)$, and hence are also on the critical line. Therefore, all roots of $\xi '(s)$ are also on the critical line. (The above argument is basically in \cite{RefPolya2}, which is in German, which we reproduce here for the convenience of readers.) 

By Theorem \ref{TheoremMeromorphic}, the attractors of  BNQN applied to $F(s)=|\xi (s)|^2/2$ are either zeros of $\xi (s)$ or zeros of $\xi '(s)$, hence belong to the critical line. 

($\Leftarrow $) First, we recall that by Theorem \ref{TheoremMeromorphic}, if $S$ is an attractor for the dynamics of BNQN applied to $F(s)$, then $S$ must be a point and moreover, must be a root of $\xi (s)\xi '(s)$. In addition, we have: 

{\bf Claim:} If $z^*$ is a root of $\xi (s)$, then $z^*$ is an attractor for the dynamics of BNQN. 

Proof of Claim: Indeed, Lemma A.2 in \cite{RefTT} shows that near $z^*$, inequality (3.1) in the paper \cite{absil-mahony-andrews} is satisfied. Therefore, since $z^*$ is an isolated (global) minimum of the function $F$, the last part of the proof of Theorem 3.2 in \cite{absil-mahony-andrews} shows that there is $\epsilon >0$ (which can be chosen so small that $z^*$ is the only root of $\xi (s)\xi '(s)$ inside the set $\{z:~|z-z^*|\leq 2\epsilon \}$) so that if $z_0\in \{z:~|z-z^*|<\epsilon \}$, then $z_n\in \{z:~|z-z^*|<2 \epsilon \}$ for all $n$. Therefore, by Theorem \ref{TheoremMeromorphic}, the sequence $\{z_n\}$ must converge to $z^*$. Hence, $z^*$ is an attractor for the dynamics of BNQN applied to $F(s)$. 

From the Claim, if all of attractors of BNQN belong to the critical line, then specially all roots of the Riemann xi function belongs to the critical line. Therefore, the Riemann hypothesis follows. 

\end{proof}

\begin{proof}[Proof of Theorem \ref{Theorem2}] 
Under the assumptions of the theorem, if $\xi '(z^*)=0$, then $z^*$ is a saddle point of $F(x,y)=|\xi (x+iy)|^2/2$. Therefore, we can apply Theorem \ref{TheoremMeromorphic} to obtain Theorem \ref{Theorem2}.  
\end{proof}

\section{Experimental results} In this section we present some experimental results concerning applying BNQN to the function $F(x,y)=|\xi (x+iy)|^2/2$. 

We first mention two difficulties concerning computations with transcendental functions like the Riemann xi function. First, transcendental functions or numbers cannot be accurately calculated on computers. Therefore, one has to use approximations. Specially, for the Riemann zeta/xi function, an effective tool is provided in \cite{RefReyna}, which is implemented in the python library mpmath \cite{Refmpmath}. Second, there are many points in the critical strip which are not close to the critical line but at which the value of the Riemann xi function is extremely small.    

The above two difficulties require one to be very carefully when numerically finding roots of the Riemann xi function. In particular, one has to use a machine precision much higher than usually needed (in the experiments here we need to choose the values 100, 300 or 1000 for the mpmath parameter mp.dps), and the experiments will take a long time to run. We also have to restrict the number of iterations to about tens to one hundred. An implementation in Python of BNQN accompanies the paper \cite{RefTT}, and this paper incorporates \texttt{mpmath} to deal with the transcendental function \(\xi\) and the high machine precision. 

{\bf Experiment 1:} We compare the basin of attractions for various methods in finding the first 8 roots of the Riemann xi function. Here the roots are: 

Root 1 $\sim$ $0.5+14.13472514173i$, with points in the basin of attraction having green colour. 

Root 2 $\sim$ $0.5-14.13472514173i$, with points in the basin of attraction  having yellow colour.  

Root 3 $\sim$ $0.5+21.02203963877i$, with points in the basin of attraction having blue colour.  

Root 4 $\sim$ $0.5-21.02203963877i$, with points in the basin of attraction having red colour. 

Root 5 $\sim$ $0.5+25.01085758014i$, with points in the basin of attraction having  pink colour. 

Root 6 $\sim$ $0.5-25.01085758014i$, with points in the basin of attraction having cyan colour.

Root 7 $\sim$ $0.5+30.42487612585i$, with points in the basin of attraction having orange colour.

Root 8 $\sim$ $0.5-30.42487612585i$, with points in the basin of attraction having purple colour.

All other points in the grid have black colour. Among other options, they may represent initial points where the sequence converge to another root or to infinity or to strange attractors. By Theorem \ref{TheoremMeromorphic}, BNQN has no strange attractors, but it is unknown whether Newton's method or Random Relaxed Newton's method can have. 

In this experiment, we choose initial points in a grid of the size $250\times 250$ points. We also rescale the y-axis by a factor of $0.1$. For example, the point $0.1+1.1i$ in Figure  corresponds to the point $0.1+11i$ in the usual coordinate system.  

\begin{figure}
    \centering
    
    \includegraphics[width=2cm]{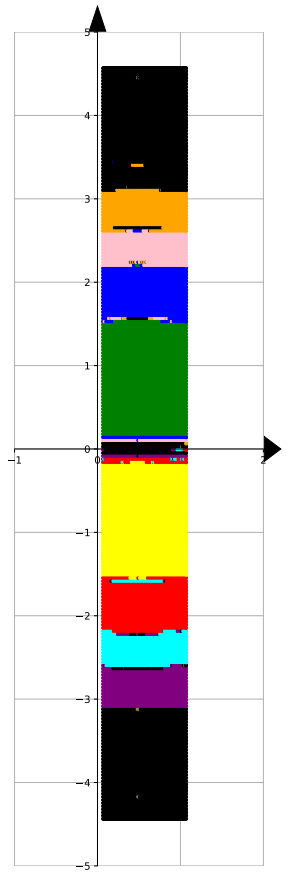}
    \includegraphics[width=2cm]{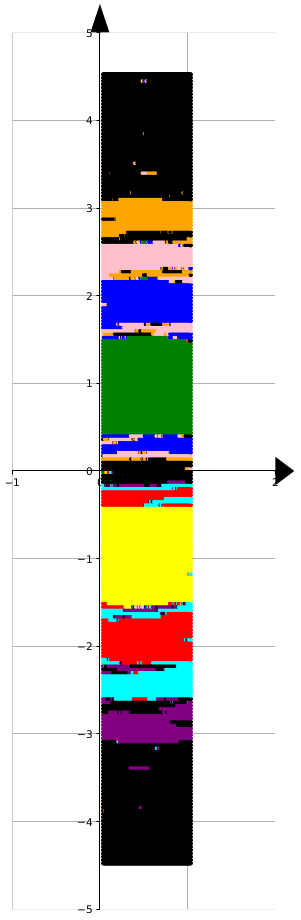}

    \bigskip
   \includegraphics[width=13cm]{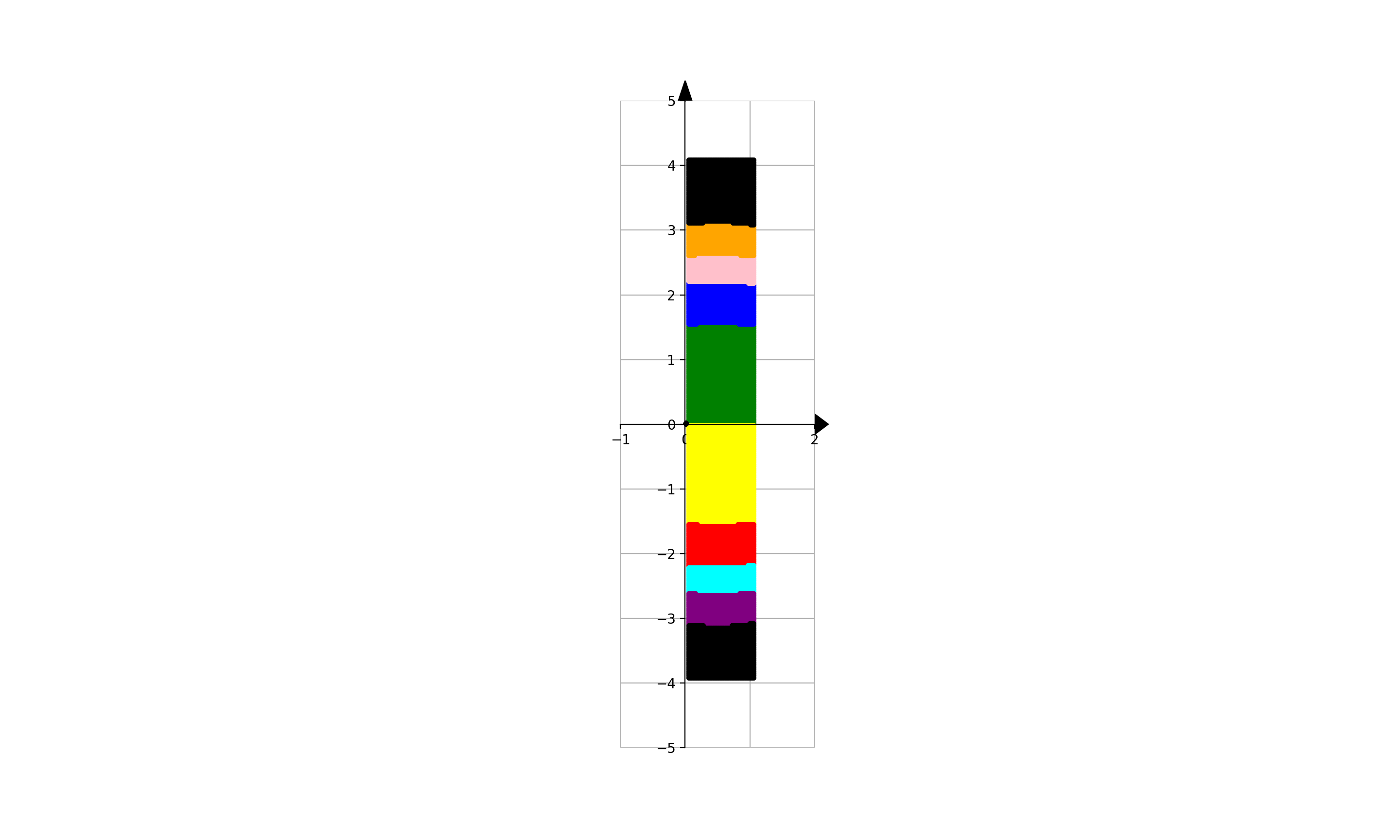}
    
    \caption{Experiment 1: Basins of attraction when applying BNQN to find roots of $\xi (z)$, specially the first 8 roots. The y-axis is rescaled by $0.1$, hence $y=1$ in the pictures correspond to $y=10$ in reality.  Pictures are referenced to from top to bottom, from left to right. Row 1: left picture is Newton's method, right picture is Random Relaxed Newton's method. Row 2: BNQN.\\ One can see that pictures for Newton's method and Random Relaxed Newton's method display fractal features, while BNQN picture is very regular and similar to Voronoi's diagrams.}
    \label{fig:f2}
\end{figure}

See Figure \ref{fig:f2} for detail. 

{\bf Experiment 2:} We test using BNQN for finding roots of the Riemann xi function, with initial points of the form $(0,10^9+(j/30))$, for $j=0,1,\ldots ,30$.

By using the function nzeros(t) in mpmath, we find that in the domain $0<\mathcal{R}(z)<1$, $10^9<\mathcal{I}(z)<10^9+1$, the Riemann zeta function has 3 roots. The function zetazero() in mpmath gives the 3 roots as follows: 
\begin{eqnarray*}
&&(0.5, 1000000000.1156508900208481613883...),\\
&&(0.5,1000000000.43402689589474340890360...),\\
&&(0.5, 1000000000.5303428567293815535254...).
\end{eqnarray*}

With initial points of the above form, BNQN can find only 1 root inside $0<\mathcal{R}(z)<1$, $10^9<\mathcal{I}(z)<10^9+1$, together with a couple of roots outside the domain. We will discus in Example 4 on how to look for the missing roots, by choosing other initial points.  

With initial point $z_0=(0, 10^9)$, after $n=30$ iterations one arrives at the point $z_{n}=(0.5, 1000000000.434026895894743408903...)$. It is checked by using sign changes of the Riemann xi function on the critical line that there is a root of the Riemann zeta function on the critical line within a $10^{-6}$ distance to $z_{n}$. 

With initial points $z_0=(0, 10^9+j/30)$, where $1\leq j\leq 24$, after $n=30$ iterations one arrives at the point $z_{n}$ which is close to that for the initial point $(0,10^9)$. 

With initial point $z_0=(0, 10^9+25/30)$, after $n=30$ iterations one arrives at the point $$z_{n}=(0.499999997983191137...,1000000001.2937395161468348297...).$$  
It is checked by using sign changes of the Riemann xi function on the critical line that there is a root of the Riemann zeta function on the critical line within a $10^{-6}$ distance to $z_{n}$. 

With initial point $z_0=(0, 10^9+26/30)$, after $n=30$ iterations one arrives at the point $$z_{n}=(0.5,1000000001.29373949270440769...).$$  
This point is within $10^{-6}$ distance to the one for the initial point $(0, 10^9+25/30)$, so corresponds to the same root. 

With initial point $z_0=(0, 10^9+27/30)$, after $n=30$ iterations one arrives at the point $$z_{n}=(0.5,1000000001.29373949270440769...).$$  
This point is within $10^{-6}$ distance to the one for the initial point $(0, 10^9+25/30)$, so corresponds to the same root. 

With initial point $z_0=(0, 10^9+28/30)$, after $n=30$ iterations one arrives at the point $$z_{n}=(0.5,1000000001.9098505173208544921453...).$$  
It is checked by using sign changes of the Riemann xi function on the critical line that there is a root of the Riemann zeta function on the critical line within a $10^{-6}$ distance to $z_{n}$. 

With initial point $z_0=(0, 10^9+29/30)$, after $n=30$ iterations one arrives at the point $$z_{n}=(0.5,1000000001.293739492704407697692...).$$  
This is close to that for the initial point $(0, 10^9+27/30)$. 

With initial point $z_0=(0, 10^9+30/30)$, after $n=30$ iterations one arrives at the point $$z_{n}=(0.5,1000000001.6054652819648720667181...).$$  
It is checked by using sign changes of the Riemann xi function on the critical line that there is a root of the Riemann zeta function on the critical line within a $10^{-6}$ distance to $z_{n}$. 

{\bf Experiment 3:} We test using BNQN for finding roots of the Riemann xi function, with initial points of the form $(0,10^{10}+(j/30))$, for $j=0,1,\ldots ,30$. 

By using the function nzeros(t) in mpmath, we find that in the domain $0<\mathcal{R}(z)<1$, $10^{10}<\mathcal{I}(z)<10^{10}+1$, the Riemann zeta function has 3 roots. The function zetazero() in mpmath gives the 3 roots as follows: 

\begin{eqnarray*}
&&(0.5,10000000000.0606343467918523955294...),\\
&&(0.5,10000000000.2802888360631074507330...),\\
&&(0.5,10000000000.7065048231449404100685...).
\end{eqnarray*}

In this case, with the above initial points, BNQN can find all the above 3 roots, together with a couple of other roots. 

With initial point $z_0=(0, 10^{10})$, after $n=30$ iterations one arrives at the point $$z_{n}=(0.4999999999999999..., 10000000000.06063434679185239552...)$$ It is checked by using sign changes of the Riemann xi function on the critical line that there is a root of the Riemann zeta function on the critical line within a $10^{-6}$ distance to $z_{n}$. 

With initial points $z_0=(0, 10^{10}+j/30)$, where $1\leq j\leq 17$, after $n=30$ iterations one arrives at the point $z_n$ which is close to that for the initial point $(0,10^{10})$.

With initial point $z_0=(0, 10^{10}+18/30)$, after $n=30$ iterations one arrives at the point $$z_{n}=(0.4999999815249915967725..., 10000000000.280288834966082790...)$$ It is checked by using sign changes of the Riemann xi function on the critical line that there is a root of the Riemann zeta function on the critical line within a $10^{-6}$ distance to $z_{n}$. 

With initial point $z_0=(0, 10^{10}+19/30)$, after $n=30$ iterations one arrives at the point $$z_{n}=(0.5, 10000000001.04055844636522984441486672855...)$$ It is checked by using sign changes of the Riemann xi function on the critical line that there is a root of the Riemann zeta function on the critical line within a $10^{-6}$ distance to $z_{n}$. 

With initial point $z_0=(0, 10^{10}+20/30)$, after $n=30$ iterations one arrives at the point $$z_{n}=(0.5, 10000000000.70650482314494041006858899...)$$ It is checked by using sign changes of the Riemann xi function on the critical line that there is a root of the Riemann zeta function on the critical line within a $10^{-6}$ distance to $z_{n}$. 

With initial point $z_0=(0, 10^{10}+21/30)$, after $n=30$ iterations one arrives at the point $z_n$, which is close to that for the initial point $(0, 10^{10}+20/30)$. 

With initial point $z_0=(0, 10^{10}+22/30)$, after $n=30$ iterations one arrives at the point $z_n$, which is close to that for the initial point $(0, 10^{10})$.

With initial point $z_0=(0, 10^{10}+23/30)$, after $n=30$ iterations one arrives at the point $$z_{n}=(0.5, 10000000001.0405584463652298444...)$$. 

With initial point $z_0=(0, 10^{10}+24/30)$, after $n=30$ iterations one arrives at the point $$z_{n}=(0.4999999999999..., 10000000001.290870692585983187939...)$$ It is checked by using sign changes of the Riemann xi function on the critical line that there is a root of the Riemann zeta function on the critical line within a $10^{-6}$ distance to $z_{n}$. 

With initial point $z_0=(0, 10^{10}+25/30)$, after $n=30$ iterations one arrives at the point $z_n$, which is close to that for the initial point $(0,10^{10}+24/30)$.

With initial point $z_0=(0, 10^{10}+26/30)$, after $n=30$ iterations one arrives at the point $z_n$, which is close to that for the initial point $(0,10^{10}+24/30)$.

With initial point $z_0=(0, 10^{10}+27/30)$, after $n=30$ iterations one arrives at the point $z_n$, which is close to that for the initial point $(0,10^{10})$.

With initial point $z_0=(0, 10^{10}+28/30)$, after $n=30$ iterations one arrives at the point $z_n$, which is close to that for the initial point $(0,10^{10}+19/30)$.

With initial point $z_0=(0, 10^{10}+29/30)$, after $n=30$ iterations one arrives at the point $z_n$, which is close to that for the initial point $(0,10^{10}+19/30)$.

With initial point $z_0=(0, 10^{10}+30/30)$, after $n=30$ iterations one arrives at the point $z_n$, which is close to that for the initial point $(0,10^{10}+19/30)$.

\textbf{Example 4:} Here we seek other initial points of the form $(0,y_0)$, for which BNQN can find the 2 missing zeros mentioned in Example 2.  

We see that for initial points of the form $(0, 10^9+j/30)$, where $j=0,\ldots ,24$, then BNQN finds the zero 
$$(0.5,1000000000.43402689589474340890360...).$$ 
On the other hand, for initial points of the form $(0,10^9+j/30)$ for $j=25,\ldots ,30$, then BNQN finds zeros in the domain $10^9+1<\mathcal{I}(z)$. 

\underline{To find the zero $(0.5, 1000000000.5303428567293815535254...)$}: The above observation suggests that we should choose some initial points of the form $(0,y_0)$ where $10^9+24/30 <y_0 <10^9+25/30$. 

We find that with two initial points $(0, 10^9+(24/30)+(1/300))$ and $(0,10^9+(24/30)+(2/300))$, then after 30 iterations BNQN is close to the zero  $(0.5, 1000000000.5303428567293815535254...)$, as wanted.

\underline{To find the zero $(0.5, 1000000000.1156508900208481613883...)$}: The above observation suggests that we should choose some initial points of the form $(0,y_0)$, where $y_0<10^9$. 

We will first look at initial points of the form $(0,10^9-j/30)$, where $ j=1,2,\ldots $, to locate an interesting position. 

For initial points of the form $(0, 10^9-j/30)$, for $1\leq j\leq 15$, after $n=30$ iterations BNQN arrives at a point $z_n$ close to the zero $$(0.5,1000000000.43402689589474340890360...).$$ 

With th initial point $(0,10^9-16/30)$, after $n=30$ iterations, BNQN is close the zero $(0.5, 1000000000.1156508900208481613883...)$, as wanted. 

\section{Some concrete ideas on using BNQN towards the Riemann hypothesis}

Here we analyse in detail the insights obtained from experiments presented in the previous section, and present some concrete ideas on using BNQN for the Riemann hypothesis. 

\subsection{Some observations}

We first report some observations from the experiments. 

{\bf Observation 1:} It seems that for all initial points $z_0$ in the critical strip, the orbit under BNQN does not converge to $\infty$. 
 
{\bf Observation 2:} It seems that BNQN can find all roots of the Rieman xi function by starting with initial points on the y-axis. 

{\bf Observation 3:} It seems that one needs at most $O(log (T))$ iterations by BNQN, from a starting point $(x_0,y_0)$ in the critical strip with $|y_0|\leq T$, to arrive in a given $\epsilon$-neighbourhood of the critical line.  

Observation 3 may be related to classical estimates for the Riemann zeta/xi function in the critical strip.

\subsection{Some concrete ideas on using BNQN for the Riemann hypothesis} We first present a precise version of the similarity between Voronoi's diagrams and basins of attraction of BNQN applied to $\xi$, based on observations from experiments. 

{\bf Conjecture 1:} Let's consider the Voronoi's diagram for all zeros of the Riemann xi function. If the intersection between  Voronoi's cell of a zero $z^*$ and the y-axis contains a non-empty open interval, then for BNQN applied to the Riemann xi function, the intersection between the y-axis and the basin of attraction of $z^*$ also contains a non-empty open interval.  

Here is a relevant consequence. 

\begin{cor} Assume that Conjecture 1 holds. Let $\mathcal{A}$ be a dense subset of the y-axis. 

- If the Riemann hypothesis holds, then BNQN applied to $\xi$ with initial points in $\mathcal{A}$ can find all zeros of $\xi$. 

- If the Riemann hypothesis fails, then BNQN applied to $\xi$ with initial points in $\mathcal{A}$ can find one zero of $\xi$ which is off the critical line. 
\label{Corollary1}\end{cor}
\begin{proof}
If the Riemann hypothesis holds, then the Voronoi's cell of any zero $z^*$ of $\xi$ will intersect the y-axis in a non-empty open interval.

If the Riemann hypothesis fails, then since the distance between any point on the critical line to the y-axis is $0.5$, while there is a zero of $\xi$ whose distance to the y-axis is $<0.5$. Therefore, there is an off-critical line zero $z^*$ of $\xi$ whose Voronoi's cell intersects the y-axis in a non-empty open interval. 

In both cases, Conjecture 1 says that starting from initial points in $\mathcal{A}$, the algorithm BNQN can find $z^*$.

\end{proof}

Here are some concrete ideas on using BNQN for the Riemann hypothesis. 

Step 1: show that Conjecture 1 holds. Of relevance, in \cite{RefFHTW} it is shown that BNQN applied to polynomials of degree $2$ produces pictures which are precisely Voronoi's diagram of the roots of the polynomial. In the proof in \cite{RefFHTW}, Armijo's condition plays an important role. 

Step 2: Show that if $z_0=(x_0,y_0)$ is in the critical strip and $\xi (z_0)\xi '(z_0)\not= 0$, then the point $z_1=(x_1,y_1)$ obtained from $z_0$ by applying BNQN to $\xi$ also satisfies $\xi (z_1)\xi '(z_1)\not= 0$ (unless $x_1=0.5$). Again, the ideas in \cite{RefFHTW} can be relevant.  As in the current literature, one can concretely try to show that $\mathcal{R}(\xi '/\xi )(z_1)\not= 0$. The first progress forward is to establish this claim for $z_0=(0,y_0)$.  An advantage with choosing an initial point $z_0=(0,y_0)$ is that we know that $\xi (z_0)\xi '(z_0)\not= 0$ (see \cite{RefCon2}).   

Step 3: Show that with an initial point $z_0=(0,y_0)$, after $n=O(\log (|y_0|))$ iterations of BNQN, one arrives at a point $z_n$ where the conditions by Kantorovich or Smale \cite{smale} for Newton's method are satisfied. The first attempt can be to establish that $|z_n-z_0|=O(\log (|y_0|))$ for all n. This is observed in the experiments reported in Section 3. In fact, it seems also that $|x_n|=O(1)$.

Step 4: Show that at a point $z_n$ which is  on the orbit of $z_0=(0,y_0)$ and which satisfies Kantorovich or Smale's conditions, then the behaviour by BNQN is the same as that for Newton's method, and the constructed sequence converges to the critical line. One idea is to consider the decreasing property for the sequence $\{|x_j-1/2|\}$, for $j\geq n$. 

The estimates needed in these steps are more flexible than those needed in Turing's approach. For example, we conjecture that if there is a method which can show that $100\%$ of the roots of the Riemann xi function lies on the critical line, then the estimates in that method can be used to make our approach successful. 

\subsection{A further idea by combining with de Bruijn- Newman's constant}  Here, we extend the ideas in the previous subsection by combining de Bruijn - Newman's constant (recalled in Section 2). From \cite{RefBruijn}\cite{RefNewman}, we know that the Riemann hypothesis is equivalent to the following: For each $t>0$, the function $H_t(z)$ has only real zeros. We can show, similarly to the case of Riemann xi function, the following: 

{\bf Claim:} For each $t>0$, the following two statements are equivalent: 1. $H_t(z)$ has only real zeros. 2. All attractors of BNQN, applied to $|H_t(z)|^2$, lie on the real axis. 

One issue when working with the Riemann xi function (or equivalently, $H_0(z)$) is that it is difficult to find a simple and good estimate for it in large domains. The situation for $H_t(z)$, for each $t>0$, is better:  there is a simple function $B_t(z)$ which asymptotically approximates $H_t(z)$, in the domain $0\leq \mathcal{R}(z)\leq 1$ (note that, the role of the real and imaginary parts in $H_t(z)$ is opposite to that in the $\xi $ function), see the introduction section in \cite{RefPolymath}. This gives hope that we will be able to proceed better when applying BNQN to the functions $H_t(z)$ with any fixed $t>0$.


\begin{thebibliography}{}

\bibitem{absil-mahony-andrews} P.-A. Absil, R. Mahony and B. Andrews, Convergence of the iterates of descent methods for analytic cost functions, SIAM J. Optim. 16 (2005), vol 16, no 2, 531--547. 

\bibitem{RefAr} L. Armijo, {\it Minimization of functions having Lipschitz continuous first partial derivatives}, Pacific J. Math. 16 (1966), no. 1, 1--3. 

\bibitem{balanzario-ortiz} E. P. Balanzario and J. S\'anchez-Ortiz, Zeros of the Davenport-Heilbronn counterexample, Mathematics of Computation 76 (2007), number 260, 2045--2049.   

\bibitem{RefBer} W. Bergweiler,  Iteration of meromorphic functions, Bulletin of the American Mathematical Society, vol 29, number 2, October 1993, pp. 151--188  (1993)

\bibitem{RefBruijn} Bruijn, de, N. G. (1950). The roots of trigonometric integrals. Duke Mathematical Journal, 17(3), 197-226.
https://doi.org/10.1215/S0012-7094-50-01720-0

\bibitem{RefBroughan1} K. Broughan, Equivalents of the Riemann Hypothesis. Volume 1: Arithmetic Equivalents. Cambridge University Press (2017)

\bibitem{RefBroughan2} K. Broughan, Equivalents of the Riemann Hypothesis. Volume Two: Analytic Equivalents. Cambridge University Press (2017)

\bibitem{RefBroughan3} K. Broughan, Equivalents of the Riemann Hypothesis. Volume 3: Further Steps towards Resolving the Riemann Hypothesis. Cambridge University Press (2023)

\bibitem{RefBui} H. M. Bui, J. B. Conrey, and M. P. Young, At least $41\%$ of the zeros of the zeta function are on the critical line, Acta Arith. 150, 35--64, 2011. 
  
\bibitem{RefCon2} J. B. Conrey,  Zeros of derivatives of Riemann's $\xi $ function on the critical line, Journal of Number Theory 16, pp. 49--74, 1983. 

\bibitem{RefCon} J. B. Conrey,  More than two fifths of the zeros of the Riemann zeta function are on the critical line, J. Reine Angew. Math. 1989 (399), 1--26. 

\bibitem{davenport-heilbronn} H. Davenport and H. Heilbronn, On the zeros of certain Dirichlet series, I and II, J. London Math. Soc. 11 (1936), 181--185, 307--312.  

\bibitem{fornaess-etal2} J. E. Fornæss, M. Hu, T. T. Truong and T. Watanabe, Backtracking New Q-Newton's method, Newton's flow, Voronoi's diagram and Stochastic root finding, arXiv:2401.01393. Accepted in Complex Analysis and Operator Theory. 

\bibitem{RefFHTW}  J. E. Forn\ae ss, M. Hu, T. T. Truong, T. Watanabe, 
Backtracking New Q-Newton's method, Schr\"oder's theorem, and Linear conjugacy, arXiv:2312.12166. 

\bibitem{RefHinkkanen} A. Hinkkanen, On functions of bounded type, Complex Variables, vol 34, pp. 119--139, 1997.   

\bibitem{RefIvic} A. Ivic, The Riemann zeta function, John Wiley and Sons, 1985. (Reprinted by Dover in 2003).  

\bibitem{RefHK} F. von Haeseler, H, Kriete, The relaxed Newton's method for rational functions, Random Comput. Dynam., 3 (1995), 71--92.

\bibitem{RefKawahira} T. Kawahira, The Riemann hypothesis and holomorphic index in complex dynamics, Experimental Mathematics 27, no 1 (2018), 37--46.  

\bibitem{RefLevinson} N. Levinson, Zeros of derivative of Riemann $\xi$ function, Bulletin of the American Mathematical Society, vol 80, number 5, pp. 951--954, September 1974.

\bibitem{RefMe} H.-G. Meier, The relaxed Newton-iteration for rational functions: the limiting case, Complex Variables Theory Appl., 16 (1991), 239--260.

\bibitem{RefMontVaug}  
H. L. Montgomery  and R. C. Vaughan,  Multiplicative number theory. I. Classical theory. Vol. 97. Cambridge Studies in Advanced Mathematics.
Cambridge University Press, Cambridge, 2006, pp. xviii+552.

\bibitem{Refmpmath}  mpmath: a {P}ython library for arbitrary-precision floating-point arithmetic (version 1.3.0), http://mpmath.org/

\bibitem{RefNewman} C. M. Newman, Fourier Transforms with Only Real Zeros,  Proceedings of the American Mathematical Society, Vol. 61, No. 2 (Dec., 1976), pp. 245-251, American Mathematical Society, https://doi.org/10.2307/2041319

\bibitem{RefNeuberger2} J. W. Neuberger, C. Feiler, H. Maier and W. P. Schleich, The Riemann hypothesis illuminated by the Newton flow of $\zeta$, Invited comment, Physica Scripta, volume 90, 108015, 2015. 

\bibitem{RefNeuberger1} J. W. Neuberger, C. Feiler, H. Maier and W. P. Schleich, Newton flow of the Riemann zeta function: separatrices control the appearance of zeros, New Journal of Physics, volume 16, 103023, 2014. 

\bibitem{RefN} Wikipedia page on Newton's fractal, https://en.wikipedia.org/wiki/Newton$\_$fractal

\bibitem{RefOdlyzko} A. Odlyzko, The $10^{21}$-st zero of the Riemann zeta function, note for the informal proceedings of the September 1998 conference on the zeta function at the Edwin Schroedinger institute in Vienna.  https://www-users.cse.umn.edu/$\sim$odlyzko/unpublished/zeta.10to21.pdf

\bibitem{RefPlattTrudgian} D. Platt, T. Trudgian, The Riemann hypothesis is true up to \(3 \cdot 10^{12}\), Bulletin of the London Mathematical
Society, Volume 53, issue 3, 792–797, 2021.


\bibitem{RefPolya2} G. Pólya, Bemerkung zur Theorie der ganzen Funktionen,  Deutsch. Math.-Verein. 24, 392--400, 1915. 

\bibitem{RefPolya} G. Pólya, Über trigonometrische Integrale mit nur reellen Nullstellen, Journal für die reine und angewandte Matematik, vol. 158(1927), pp. 98-99.

\bibitem{RefPolymath} D. H. J. Polymath, Effective approximation of heat flow evolution of the Riemann \(\xi\) function, and a new upper bound for the de Bruijn-Newman constant, Research in the mathematical sciences, 6 (3), paper no 31, 67 pp, 2019.

\bibitem{RefPratt} K. Pratt, N. Robles, A. Zaharescu and D. Zeindler, More than five twelfths of the zeros of the Riemann zeta function are on the critical line, Research in the Mathematical science, volume 7, article 2, 2020. 


\bibitem{RefReyna} J. Arias de Reyna, High precision computation of {R}iemann's zeta function by the {R}iemann-{S}iegel formula, {I}, Mathematics of Computation, vol 80, number 274, 995--1009, 2011. 

\bibitem{RefRodgersTao} B. Rodgers and T. Tao, The de Bruijin-Newman constant is non-negative, Forum of Mathematics, 8, e6.

\bibitem{RefSch} W. P. Schleich, I. Bezdekova, M. B. Kim, P. C. Abbott, H. Maier, H. L. Montgomery and J. W. Neuberger, Equivalent formulations of the Riemann hypothesis based on lines of constant phase, Physica Scripta, 93, 065201 (11 pp), 2018.

\bibitem{RefSchr} D. Schleicher, Newton's method as a dynamical system: efficient root finding of polynomials and the Riemann zeta function, Fields Inst. Commun. 53 (2008), 213--224.

\bibitem{spira} R. Spira, Some zeros of the Titchsmarsh counterexample, Math. Comp. 63 (1994), no 208, 747--748.  

\bibitem{RefS} H. Sumi, Negativity of Lyapunov exponents and convergence of generic random polynomial dynamical systems and random relaxed Newton's method, Communications in Mathematical Physics, vol 384, 1513--1583, 2021. 
      
\bibitem{RefT} T. T. Truong, Backtracking New Q-Newton's Method: A Good Algorithm for Optimization and Solving Systems of Equations (2023)	arXiv:2209.05378
		
\bibitem{RefTT} T. T. Truong, T. D.  To, H.-T. Nguyen, T.  H. Nguyen, H.  P. Nguyen, M. Helmy, A fast and simple modification of Newton's method avoiding saddle points. J Optim Theory Appl (2023). https://doi.org/10.1007/s10957-023-02270-9

\bibitem{smale} S. Smale, {\it Newton's method estimates from data at one point}, in The merging of disciplines: new directions in pure, applied and computational mathematics (Laramie, Wyoming, 1985), pp. 185--196, Springer, New York, 1986. 

\bibitem{RefV1} G. Voronoi, Nouvelles applications des parametres continus a la theorie des formes quadratiques, Premier memoire, Sur quelques proprietes des formes quadratiques positive parfaites, Journal fur die Reine und Angwandte Mathematik 1908 (133), pp. 97--178, 1908. 

\bibitem{RefV2} G. Voronoi, Nouvelles applications des parametres continus a la theorie des formes quadratiques, Deuxieme memoire, Recherches sur les paralleloedres primitifs, Journal fur die Reine und Angwandte Mathematik 1908 (134), pp. 198--178, 287. 
		
\end{thebibliography}
\end{document}